\documentclass[12pt]{article}      

\usepackage[english]{babel}
\usepackage{amssymb,amsmath,latexsym}

\newtheorem{thm}{Theorem}
\newtheorem{lem}{Lemma}

\newtheorem{corol}{Corollary}

\title{Finite and countable infinite products of Probabilistic Normed
Spaces}

\bigskip
\author{Bernardo Lafuerza-Guill\'{e}n\\ \\
Departamento de Estad\'{\i}stica y Matem\'{a}tica
Aplicada\\Universidad de Almer\'{\i}a
        \\04120-Almer\'{\i}a, Spain}     
\date{}      

\begin{document}             
\maketitle                 
\begin {abstract}

In this work we first give for PN spaces results parallel to those
 obtained by Egbert for the product of PM spaces, and
generalize  the results by Alsina and Schweizer in \cite{alsina2}
in order to study non-trivial products and the product of $m$--transforms of several PN spaces.\\
 In addition we present a detailed study of $\alpha$--simple product PN spaces and, finally,
 the product topologies in PN spaces which are products of countable families
of PN spaces.

\bigskip


\bigskip
 A.M.S. CLASSIFICATION: 54E70
\end{abstract}
\section{\bf INTRODUCTION}
We assume that the reader is acquainted with the basic notions of
the theory of PN spaces. These, as well as terms and concepts not
defined in the body of this paper, may be found in
\cite{alsina1,alsina2,alsina3,LRS1,LRS2}.
\bigskip

\noindent DEFINITION 1. A probabilistic metric space (henceforth
and briefly, a PM space) is a triple $(S,\mathcal{F},\tau)$ where
$S$ is a nonempty set (whose elements are the points of the space
), $\mathcal{F}$ is a function from $S \times S$ into $\Delta^+$,
$\tau$ is a triangle function, and the following conditions are
satisfied for all $p, q, r$ in $S$:
\begin{itemize}
\item $\mathcal{F}(p, p) = \epsilon_0$
\item $\mathcal{F}(p, q)\neq \epsilon_0$ if $p \neq q$
\item $\mathcal{F}(p, q)=\mathcal{F}(q, p)$
\item $\mathcal{F}(p, r)\geq \tau(\mathcal{F}(p, q), \,\mathcal{F}(p,
r))$.\\

\noindent If $(S,\mathcal{F},\tau)$ is a PM space, then we also
say that $(S,\mathcal{F})$ is a PM space under $\tau$.
\end{itemize}
 \noindent DEFINITION 2. A probabilistic Normed Space,
briefly a PN space, is a quadruple $(V,\nu,\tau,\tau^*)$ in which
$V$ is a linear space, $\tau$ and $\tau^*$ are continuous
triangle functions with $\tau\leq\tau^*$ and $\nu$, the
probabilistic norm, is a map $\nu: V\to\Delta^{+}$ such that

\bigskip
\noindent
 (N1)\, $\nu_p=\epsilon_0$ if, and only if, $p=\theta$, $\theta$
being the null vector in $V$;\\ (N2)\, $\nu_{-p}=\nu_p$\quad for
every $p\in V$;\\ (N3)\, $\nu_{p+q}\geq\tau(\nu_p,\nu_q)$ for all
$p,q \in V$;\\ (N4)\, $\nu_p\leq\tau^*(\nu_{\alpha
p},\nu_{(1-\alpha)p})$ for every $\alpha\in [0,1]$ and for every
$p\in V$.

\bigskip

 If, instead of (N1), we only have
$\nu_{\theta}=\epsilon_0$, then we shall speak of a {\it
Probabilistic Pseudo Normed Space}, briefly a PPN space. If the
inequality (N4) is replaced by the equality
$\nu_p=\tau_M(\nu_{\alpha p},\nu_{(1-\alpha)p})$, then the PN
space is called a \v{S}erstnev space and, as a consequence, a
condition stronger than (N2)\,holds, namely
\bigskip
$$\forall\lambda\neq 0\:\forall p\in V\qquad \nu_{\lambda
p}=\nu_p\left(\frac{j}{|\lambda |}\right) $$

\bigskip
Here $j$\,\, is the identity map on $\mathbf{R}$,\quad i.e.
\,$j(x):=x\, (x\in \mathbf{R})$.\\
A \v{S}erstnev space is denoted by $(V,\nu,\tau)$.

\bigskip
There is a natural topology in a PN space $(V,\nu,\tau,\tau^*)$,
called the {\it strong topology}; it is defined, \,for $t>0$, by
the neighbourhoods $$\mathcal{I}_p(t):=\{q\in
V:d_S(\nu_{q-p},\epsilon_0)<t\}=\{q\in V:\nu_{q-p}(t)>1-t\}.$$

\bigskip

\noindent By setting $F\leq G$ whenever $F(x)\leq G(x)$ for every
$x\in {\bf R}^+$ and $F,G \in \Delta^+$, one introduces a
natural ordering in $\Delta^+$.\\

\noindent DEFINITION 3. Let $(V,\,\parallel\cdot\parallel )$ be a
normed space and let $G\in \Delta^+$ be different from
$\epsilon_0$ and $\epsilon_{+\infty}$; define $\nu: V \rightarrow
\Delta^+$ by $\nu_{\theta}=\epsilon_0$ and
$$
\nu_p(t):= G\left( \frac{t}{{\parallel p
\parallel^{\alpha}}}\right)
$$
where $\alpha > 0$ and $\alpha \neq 1$. Then the pair $(V,\, \nu)$
will be called the $\alpha$--simple space generated by $(V,
\,\parallel \cdot
\parallel)$ and by $G$.\\

 \noindent DEFINITION 4. Let $\tau_1,\tau_2$ be two triangle
functions. Then $\tau_1$ dominates $\tau_2$, and we write
$\tau_1\gg\tau_2$, if for all $F_1,F_2,G_1,G_2 \in \Delta^+$, $$
\tau_1(\tau_2(F_1,G_1),\tau_2(F_2,G_2))\geq
\tau_2(\tau_1(F_1,F_2),\tau_1(G_1,G_2))$$ Notice that since
  $\tau_1$ is associative one has   $\tau_1\gg\tau_1$,\, so that
\lq\lq dominates\rq\rq \,is reflexive but its transitivity is
still an open question.\\

 \noindent DEFINITION 5. Given a left-continuous $t$--norm $T$, i. e. a left-continuous
 binary operation on  [0, 1] that is commutative,\, nondecreasing in
 each variable and has 1 as identity, a triangle function\, ${\bf T}$ is
 the function defined via

 $${\bf T}(F, G)(x):= T(F(x), G(x)).$$

\bigskip

\noindent DEFINITION 6. For $b\in \left]0, +\infty \right]$, let
$M_b$ be the
 set of all continuous and strictly increasing functions $m$ from
 $[0,b]$ onto $R^+=[0,\infty]$. For any $F$ in $\Delta^+$ and any
 $m$ in $M_b$ let $Fm$ be the functions on $\mathbb{R}^+$ given by

 $$(Fm)(x)=\left\{\begin{array}{lcl}F(m(x)), &{\rm }& x\in[0,b[,\\
 \lim_{x\rightarrow b^-}F(m(x)), & {\rm  }& x=b,\\
 1, & {\rm  }& x>b,\end{array}\right.$$
 if $b\in \left]0,+\infty\right[$, and by (Fm)(x)=F(m(x)) for all $x \geq 0$
 if $b=+\infty$.\\
 If $b\in]0,+\infty[$, then $Fm$ is in $\mathcal {D}^+ $ and
 $Fm\geq\epsilon_b$. If $F$ itself is in $\mathcal {D}^+ $, then
 $Fm$ is continuous at $b$. Also, $\varepsilon_t m
 =\epsilon_{m^{-1}(t)}$
 for any $t$ in $\mathbb{R}^+$.\\

 \noindent DEFINITION 7. Let $\tau$ be a triangle function and let
 $m$ belong to $M_b$. Then $m$ is $\tau$--superadditive if, for all
 $F,G$ in $\Delta^+$,
 $$\tau(F,G)m\geq\tau(Fm,Gm).$$
 The funtion $Fm$ is called the $m$--transform of $F$.
\bigskip

 \noindent The following results Lemma 1, 2, 3 and Theorem 1 can be
seen in \cite{alsina2}.
\begin{lem}\hspace{-0.23cm}\textbf{.} If $\tau$ is a triangle function such that for all
$s,t$ in $\mathbb{R}^+$,
$$\tau(\epsilon_s,\epsilon_t)=\epsilon_{s+t},$$
and if $m$ in $M_b$ is $\tau-$ {\it superadditive}, then $m$ is
superadditive, i.e., for all $x,y\in[0,b]$,
$$m(x+y)\geq m(x)+ m(y) .$$
\end{lem}

\begin{thm}\hspace{-0.22cm}\textbf{.}
Let $T$ be a continuous t-norm and let $m\in M_b$. Then $m$ is
$\tau_T$--superadditive if and only if $m$ is superadditive (see
\cite{alsina2}, Theorem 3).
\end{thm}
\begin{lem}\hspace{-0.22cm}\textbf{.} The t-norm $W$ satisfies the following
relationship:
$$
W \left(\sum_{i=1}^{\infty}\frac{a_i}{2^i},\,
\sum_{i=1}^{\infty}\frac{b_i}{2^i} \right)\leq
\sum_{i=1}^{\infty}\frac{1}{2^i}\, W(a_i, b_i),
$$
for every $(a_i), (b_i)\in [0, 1]$.
\bigskip
\end{lem}
\begin{lem}\hspace{-0.22cm}\textbf{.} If $a_i,\, b_i\in [0, 1]$ for every $i\in \mathbf{N}$,
the inequality
$$
W^*\left(\sum_{i=1}^{\infty}\frac{a_i}{2^i},\,
\sum_{i=1}^{\infty}\frac{b_i}{2^i}\right)\geq
\sum_{i=1}^{\infty}\frac{1}{2^i}\,W^*(a_i, b_i)
$$
where $ W^*(x, y):= M \{x + y, 1\} $, holds.
\end{lem}

\section {\bf Finite $\tau$--Products of PN spaces}
\bigskip
In this section, we give our definition of a $\tau$--product of
two probabilistic normed spaces which is a generalization of a
parallel result from Egbert about the $\tau$--product of PM
spaces). Moreover we show a necessary and sufficient condition
for the $\tau$--product of two general PN spaces of \v{S}erstnev
to be a \v{S}erstnev space as well as a sufficient condition for
the $\tau$--product of two Menger PN spaces
to be also a PN space of Menger.\\
The proof of most theorems is omitted
 since it is just a matter of straightforward verification of the
 statements.
\bigskip

\noindent DEFINITION 8. Let $(V_1,\nu_1,\tau,\tau^*)$ and
$(V_2,\nu_2,\tau,\tau^*)$ be two PN spaces under the same triangle
functions $\tau$ and $\tau^*$. Let $\tau_1$ be a triangle
function. Their $\tau_1$--product is the quadruple

$$ (V_1 \times V_2,\nu^{\tau_1},\tau,\tau^*) $$ where
$$\nu^{\tau_1}:V_1 \times V_2\longrightarrow \Delta^+ $$ is a
probabilistic seminorm defined by $$ \nu^{\tau_1}(p,q):=
\tau_1(\nu_1(p),\nu_2(q)) $$ for any $(p,q)\in V_1 \times V_2$.

\begin{thm}\hspace{-0.22cm}\textbf{.} Let $(V_1,\nu_1,\tau,\tau^*)$,
$(V_2,\nu_2,\tau,\tau^*)$ and $\tau_1$ be two PN spaces under the
same triangle functions and a triangle function $\tau_1$
respectively. Assume that $\tau^* \gg \tau_1$ and $\tau_1 \gg
\tau$, then the $\tau_1$--product $(V_1 \times
V_2,\nu^{\tau_1})\,$is a PN space under $\tau$ and $\tau^*$.
\end{thm}

\bigskip
\noindent
 {\bf Example 1}. The $\bf T$--product $(V_1 \times
V_2,\,\nu^ {\bf T})$ \, of\, the\,  two PN spaces $(V_1,\,\nu_1,
\tau_T, \,\bf M)$ and $(V_2,\,\nu_2,\tau_T,\, \bf M)$ is a PN
space under $\tau_T$ and $\bf M$.

\bigskip
 \noindent
 {\bf Example 2}.
Let $(V_1,F,\,\bf M)$ and $(V_2,G,\,\bf M)$ be two equilateral PN
spaces with distribution functions F, G respectively. Then, their
{\bf M}--product is an equilateral PN space with a d.d.f. given by
${\bf M}(F,\,G)$.\\ In particular, if $F \equiv G$, the {\bf
M}--product is an equilateral PN space with the same distribution
function $F$.
\bigskip

One may wonder whether the PM space associated with the
$\tau_1$--product of two PN spaces
 characterized in Theorem 2 coincides with the $\tau_1$--product of
 corresponding PM spaces. The following theorem gives an answer in the
 affirmative to this question.\\

\begin{thm}\hspace{-0.22cm}\textbf{.} Let $(V_1,\,\nu_1,\tau,\tau^*)$ and
$(V_2,\,\nu_2,\tau,\tau^*)$ be two PN spaces under the same
triangle functions $\tau$ and $\tau^*$ and let the triangle
function $\tau_1$ be such that $\tau^*\gg\tau_1$ and
$\tau_1\gg\tau$. Let $(V_1,\mathcal{F}_1,\tau)$ and
  $(V_2,\mathcal{F}_2,\tau)$ be the same spaces regarded as
 PM spaces. Then, the $\tau_1$--product $(V_1
\times V_2,\, \nu^{\tau_1},\,\tau,\, \tau^*)$ regarded as a PM
space coincides with the $\tau_1$--product $(V_1 \times
V_2,\,\mathcal{F}^{\tau_1})$.
\end{thm}

\bigskip
It is known that if $(S_1,\,\mathcal{F}_1)$ and $(S_2,\,
\mathcal{F}_2)$ are the simple spaces $(S_1, d_1, G)$ and $(S_2,
d_2, G)$, respectively,\, and if $d_{Max}$ is the metric on $S_1
\times S_2$ defined by
$$
d_{Max}((p_1, p_2), (q_1, q_2))=Max (d_1(p_1, q_1), d_2(p_2, q_2))
$$
then the {\bf M}--product of $(S_1,\,\mathcal{F}_1)$ and $(S_2,\,
\mathcal{F}_2)$ is the simple space $(S_1 \times S_2, d_{Max}, G)$
(see \cite{schw-sk1}, p. 211, Theorem 12.7.8). But in principle,
if one has two simple PN spaces $(V_1,(\|\cdot\|)_1,G,M)$ and
$(V_2,(\|\cdot\|)_2,G,M)$ with the same d.d.f. $G$ its {\bf
M}--product is not necessarily a PN space  because of the
assumption $\tau_{M^*}\gg {\bf M}$ of Theorem 1 fails here.\\
 Now then, when one replaces $\tau_{M^*}$ by $\bf{M}$ in these simple
PN spaces we obtain the PN spaces $(V_1,\|\cdot\|_1,G,{\bf M})$
and $(V_2,\|\cdot\|_2,G,{\bf M})$ respectively, as is easily
checked. Both of them are \v{S}erstnev. Can the {\bf M}--product
of these be a simple PN space?

\noindent
   The following theorem answers that question in the affirmative.

\bigskip
\begin{thm}\hspace{-0.22cm}\textbf{.} Let $(V_1,\|\cdot\|_1,G,\tau_M,{\bf M})$ and $(V_2,\|\cdot\|_2,G,\tau_M,
{\bf M})$ and $\|\cdot\|_3$ be the two above mentioned PN spaces
and the norm defined on $V_1 \times V_2$ by

$$
\|\bar {p}\|_3:=
\|p_1\|_1\vee\|p_2\|_2,
$$
 with $\bar {p}=(p_1,p_2)\in V_1 \times
V_2$. Then $(V_1 \times V_2,\|\cdot\|_3,G,\tau_M,{\bf M})$\,is a
simple PN space that coincides with the ${\bf M}$--product of the
given simple spaces. Furthermore it is \v{S}erstnev.
\end{thm}

\bigskip
Can a $\tau$--product of two simple PN spaces with the same
generator function $G$  be a simple PN space also with $G$ as
generator ?. The following theorem answers this question in the
affirmative.

\bigskip
\begin{thm}\hspace{-0.22cm}\textbf{.} The $\tau_M$--product of two simple PN spaces $(V_1,\|\cdot\|_1, \,G,\,M)$ and
$(V_2,\|\cdot\|_2,\,G,\,M)$ is the simple space under $M$
generated by $(V_1 \times V_2,\|\cdot\|_s)$ and the same d.d.f.
$G$, namely, $(V_1 \times V_2,\|\cdot\|_s,\,G,\,M)$, where
$\|\cdot\|_s$ is the classic norm defined via
$$\|\cdot\|_s:=\|\cdot\|_1 + \|\cdot\|_2. \qquad \Box$$
\end{thm}

\begin{thm}\hspace{-0.22cm}\textbf{.}  Let $(V_1,\,\nu_1,\tau)$ and
 $(V_2,\,\nu_2,\tau)$ be two \v{S}erstnev spaces under the same triangle function $\tau$. Let us
assume that $\tau_1$ is a triangle function such that
 $\tau_1\gg\tau$, then, their
$\tau_1$--product is also a \v{S}erstnev space if, and only if
$\tau_1\gg\tau_M$ and $\tau_M\gg\tau_1$.
\end{thm}

\noindent
 {\bf Proof} \hspace{-0.1cm}\textbf{:} By Theorem 1 the $\tau_1$--product of the two \v{S}erstnev spaces exists. Now since both PN spaces
 are \v{S}erstnev one has, for all $\alpha \in [0, 1]$ and for
 all
$(p,q)\in V_1\times V_2$,\\

\begin{eqnarray*}
 & \nu^{\tau_1}(p,q)=\tau_1(\nu_1(p), \nu_2(q)) \\
 &=\tau_1[\tau_M(\nu_1(\alpha p),\nu_1((1-\alpha)p)),\tau_M(\nu_2(\alpha
 q),\nu_2((1-\alpha)q))] \\
& =\tau_M [\nu^{\tau_1}(\alpha(p,q)),\nu^{\tau_1}((1-\alpha)(p,q))]\\
&= \tau_M[\tau_1(\nu_1(\alpha p),\nu_2(\alpha
 q)),\tau_1(\nu_1((1-\alpha)p),\nu_2((1-\alpha)q))];
\end{eqnarray*}
 this implies
$\tau_1\gg\tau_M$  and $\tau_M \gg \tau_1$.\\

Conversely, let $\tau_1\gg\tau_M$  and $\tau_M\gg\tau_1$;
then\\
\begin{eqnarray*}
& \tau_1(\nu_1(p), \nu_2(q))\\
&= \tau_1[\tau_M(\nu_1(\alpha p),
\nu_1((1-\alpha)p)),\tau_M(\nu_2(\alpha q), \nu_2((1-\alpha
)q))]\\
&=\tau_M[\tau_1(\nu_1(\alpha p),\nu_2(\alpha q)),
\tau_1(\nu_1((1-\alpha)p)),\nu_2((1-\alpha)q))]\\
&= \tau_M[\nu^{\tau_1}(\alpha(p, q)), \nu^{\tau_1}( (1-\alpha)(p,
q ))],
\end{eqnarray*}
hence the assertion.\quad $\Box$

\begin{corol}\hspace{-0.22cm}\textbf{.} Under the same assumptions as in Theorem 5, if $\tau \equiv \tau_T$ then the
$\tau_1$--product $(V_1\times V_2,\nu^{\tau_1},
 \tau_T)$ is a Menger space.
 \end{corol}
\bigskip

 \noindent
 Now, the $\tau_1$--product of two Menger spaces can
 again be a Menger space. A sufficient condition is provided by the
 following theorem.
 \begin{thm}\hspace{-0.22cm}\textbf{.}
 Let $(V_1,\nu_1,T)$ and  $(V_2,\nu_2,T)$ be two Menger PN spaces. If
 $T_0$ is a left-continuous t-norm   that satisfies the
 conditions\,
 $T^*\gg T_0\quad \mbox{and}\quad T_0 \gg T,$
 then the $\tau_{T_0}$--product

 $$
 (V_1\times V_2,\nu_1\tau_{T_0}\nu_2,\tau_T, \tau_{T^*})
 $$
\noindent
 is a Menger PN space under T.
\end{thm}
\noindent
 {\bf Proof}: It suffices to apply Lemma 12.7.3 in \cite{schw-sk1}
and Theorem 1. $\Box$ \\

\noindent {\bf Example 3}. It is known that for every $t$--norm
$T$ one has $M \gg T$ (this is a result due to R. Tardiff) and it
is easily checked that $T^* \gg M^*$. Then the $\tau_M$--product
of two Menger PN spaces is also a Menger PN
space.\\

In the next section we study non trivial products.

 \bigskip
 \section {\bf Countable $\tau$--Products of PN spaces}

We shall need some preliminaries before stating the main results
of this section.\\

\bigskip
Let $(V, \parallel \cdot \parallel)$ be a normed space, $G\in
\mathcal {D}^+ $ a strictly increasing continuous d.d.f., $T$ a
strict t-norm with additive generator $f$. It is known (see
\cite{LRS2}; Lemma 3.2) that when $\alpha \in ]1, +\infty[$, then
$(V,
\parallel \cdot
\parallel, G; \alpha )$ is a Menger PN space under $T$ if, and only
if, for $s,t\in ]0, +\infty[$ and for all $p, q \in V$ such that
$p\neq \theta$, $q\neq \theta$, $p + q \neq \theta$, the
inequality
$$\parallel p + q \parallel^{\alpha}(f \circ G)^{-1}(s + t)\leq \parallel p \parallel^{\alpha}
(f \circ G)^{-1}(s)+ \parallel q \parallel^{\alpha}(f \circ
G)^{-1}(t)$$ holds.\\
\noindent

Let $(V, \parallel \cdot \parallel)$ be a normed space and let
$\alpha > 1$. If the d.f. $G \in \mathcal {D}^+$ is continuous and
strictly increasing, then (\cite {LRS2}; Section 3)  $(V,
\parallel \cdot \parallel, G; \alpha)$ is a Menger PN space under
the strict t-norm defined for all $x, y$ in $[0, +\infty]$ by
$$
T_G(x, y):= G \left(\{[G^{-1}(x)]^{\frac{1}{1- \alpha}} +
[G^{-1}(y)]^{\frac{1}{1- \alpha}}]\}^{1 - \alpha}\right)
$$
\bigskip

\noindent One is now ready to state the main results of this
section. The following one is the analogue of Theorem 12.7.9 in
\cite{schw-sk1} and shows the relevance of the t-norm $ T_G $ in
countable $ \tau $-products.\\
But in order to state it one gives a previous definition.
\bigskip

 \noindent DEFINITION 9. The mapping $\parallel \cdot \parallel_{\beta} :
  V_1 \times V_2 \rightarrow \bf{R}^+$ is defined for all
 $\bar{p} =(p_1, p_2)\in V_1 \times V_2$ via

  $$
\parallel \bar{p} \parallel_{\beta}:= \left( \parallel p_1 \parallel_{1}^{\beta}
 + \parallel p_2 \parallel_{2}^{\beta}\right)^{\frac{1}{\beta}}.
 $$

\begin{thm}\hspace{-0.22cm}\textbf{.} Let $(V_1,\parallel \cdot \parallel_1,\, G;\, \alpha)$ and $(V_2,
\parallel \cdot \parallel_2,\,G;\, \alpha)$ be two Menger PN spaces
under $ T_G $ and let $\alpha > 1$. If the d.f. $G \in \mathcal
{D}^+ $ is continuous and strictly increasing, then, the ${\bf
T}_G$--product $(V_1 \times V_2,
\parallel \cdot \parallel_{\beta},\, G;\,\alpha)$ of the given PN spaces is a Menger PN space
under $ T_G $.

\end{thm}

\noindent
 {\bf Proof}: It suffices to notice that the mapping in definition
 9 is a norm for all  $ \beta\in \,\, ]0, + \infty[ $. $\Box$\\

\noindent
 The space $(V_1 \times V_2,
\parallel \cdot \parallel_{\beta},\, G;\, \alpha)$ with $\beta=\frac{\alpha}{\alpha -1}$\, is
$ \alpha $-simple:\,\,it suffices to show that\\
\begin{eqnarray*}
&  \nu^{\bf {T}_G}(\bar{p})(j)= T_{G} \left(G\left(
\frac{j}{\parallel \, p_{1} \,\parallel_{1}^{\alpha}}\right), \,G
\left( \frac{j}{\parallel \,p_{2}
 \,\parallel_{2}^{\alpha}}\right) \right)=\\
& G  \left(\frac{j}{\left(\parallel \, p_{1}\,
\parallel_{1}^{\frac{\alpha}{\alpha -1}} +
\parallel \, p_{2}\, \parallel_{2}^{\frac{\alpha}{\alpha -1}}\right)^{\alpha -1}
}\right)= G\left( \frac{j}{{\parallel \,\bar{p}
 \,\,\parallel}_{\beta}^{\alpha}}\right).
\end{eqnarray*}
\noindent
 where $j$ \,is the identity map on ${\bf R}$
\bigskip

The analogous theorem for the case $\alpha\in ]0, 1[$ is an open
problem (see \cite {LRS2}; Section 3).\\

\noindent We know that if $\alpha > 1$ there exist normed spaces,
$(V,\, \parallel\cdot\parallel)$ with the following property: ``
if $G\in \Delta^+$ is continuous and strictly increasing, then the
t-norm $T_G$ is the strongest continuous t-norm under which
$(V,\, \parallel\cdot\parallel,\, G;\, \alpha)$ is a Menger PN
space ''(see \cite {LRS2}; Theorem 3.3). However, a new
phenomenon arises in the case of product PN spaces, for contrary
to the above, in this case the $t$--norm $T_G$ is not the
strongest continuous $t$--norm under which $(V_1 \times V_2,
\parallel \cdot \parallel_{\beta},\, G;\, \alpha)$ is a Menger PN
space, as is easily checked. Such a phenomenon is today an open
problem.
\bigskip

 In the sequel we study a special kind of probabilistic norms on the countable
 product of a family of PN spaces.
 \bigskip

\begin{thm}\hspace{-0.22cm}\textbf{.} Let $(V,\nu,\,\tau_T,\, {\bf T})$ be a PN space
and suppose $m$ in $M_b$ is $\tau_T-${\it superadditive}. Let
$\nu_m$ be the map defined for any $p$ in $V$ by
$$\nu_m(p):= \nu_p m,$$
 and $\nu_p=\nu(p)$. Then the m-transform
$(V,\nu_m,\tau_T,{\bf T})$ is a PN space.
\end{thm}

\noindent {\bf Proof}:
(N1) and (N2) are evident.\\

 (N3) Since $(V,\nu,\tau_T, {\bf T})$ is a PN space and $m$ is
 superadditive, for any $p$ in $V$, one has
 $$\nu_{p+q}m\geq \tau_T(\nu_p,\nu_q)m\geq \tau_T(\nu_p m,\nu_q m).$$

 (N4)
 \begin{eqnarray*}
& \nu_p m(x)\leq {\bf T}(\nu_{\alpha p}, \nu_{{(1-\alpha)
}p})m(x)=T(\nu_{\alpha p}m(x),
  \nu_{{(1-\alpha) }p}m(x)
 )\\
 & ={\bf T}(\nu_{\alpha p}m, \nu_{{(1-\alpha) }p}m)(x),
 \end{eqnarray*}
 for every $\alpha \in[0,1]$, for every $p\in V$ and for every
 $x\in\mathbb{R}^+$, whence
 $$\nu_p m\leq {\bf T}(\nu_{\alpha p}m, \nu_{{(1-\alpha) }p}m),$$
  for every $\alpha \in[0,1]$ and for every $p\in V$.
\bigskip

\begin{corol}\hspace{-0.22cm}\textbf{.} Let\,  $T_1, T_2$ be two $t$-norms such that\, $T_1
\leq T_2$, then the $m$--transform of any one of PN spaces
$(V,\nu,\tau_{T_1},\tau_{T_2})$ is a PN space under $\tau_{T_1}$
and ${\bf {T}_2}$.
\end{corol}

\noindent {\bf Proof}: \, If \,$T_1 \leq T_2$ one has
$\tau_{T_1}\leq\tau_{T_2}$, and now it suffices to apply to the
axiom (N4) the well-known inequality $\tau_T\leq {\bf T}$ for any
$t$--norm $T$ (see \cite{schw-sk1}).
\bigskip

 Now let us recall some conventions and results about infinite
 $\tau$-products.\\
 Since the $\tau_T$ operations are associative, for any sequence
 ${F_i}\in \Delta^{+}$, the $n$--fold  $\tau_T$--product $\tau_T^n (F_1,\ldots,F_{n+1})$
is well defined for each $n$ as the serial iterates of $\tau_T$,
defined recursively via

$$ Dom \,\tau_T^n =(\Delta^+)^{(n+1)},\qquad \tau_T^1=\tau_T $$
\noindent
 and

$$
\tau_T^{n+1}(F_1,\ldots,F_{n+1},F_{n+2})=\tau_T(\tau_T^n
(F_1,\ldots,F_{n+1}),F_{n+2}).
$$
\bigskip

If $T$ is a continuous $t$--norm, then it is well known that
\begin{eqnarray*}
&\lim_{n\to+\infty}\tau_T^n(F_1,\ldots,F_{n+1})(x)\\
&=\sup\{\lim_{n\to+\infty}T^n(F_1(x_1),\ldots,F_{n+1}(x_{n+1}))\},
\end{eqnarray*}
where the supremum is taken with respect to all sequences
$\{x_n\}$
of positive numbers such that $\sum_{i=1}^{\infty}x_n=x$.\\
Furthermore $\tau_T^{\infty}$ is defined on sequences $\{F_n\}$ in
$\Delta^+$ by
$$\tau_T^{\infty}{F_n}(x)=l^-(\lim_{n\to+\infty}\tau_T^n(F_1,\ldots,F_{n+1})(x)).$$

\bigskip
\noindent DEFINITION 10. Let $\{(V_i,\nu^i,\tau_T,{\bf {T}}) |\,
i\in \mathbf{N}\}$ be a countable family of proper PN spaces,
i.e., PN spaces with a continuous triangle function $\tau$ that
satisfies $\tau(\epsilon_s,\epsilon_t)\geq\epsilon_{s+t}$. The
function $\tau_T$ is one of these. Let $b_i$ be an infinite
sequence of positive numbers  such that the series
$\sum_{i=1}^{\infty}b_i$ converges. For each $i\in \mathbf{N}$
one chooses a function $m_i \in M_{b_i}$ that is
$\tau_T$--superadditive. Let $\{(V_i,G^i,\nu^i m_i,\tau_T,{\bf T})
|\, i\in \mathbf{N}\}$ be the $m_i-$transform of the family given
before, in which the map $G^i: V_i\rightarrow \mathbb{R}^+$ is
defined by $G_{p_i}^i:=\nu_{p_i}^i m_i$. We adopt the convention
$V=\prod_{i=1}^{\infty}V_i$. Now, for any sequence
$\bar{p}=(p_i)\in \prod_{i=1}^{\infty}V_i$ for all $i\in
\mathbf{N}$, we set $G_{\bar{p}}=\tau_T^{\infty}G_{p_i}^i$.

\begin{lem}\hspace{-0.22cm}\textbf{.} The function $G_{\bar{p}}$ defined as in Definition 7 is in
$\mathcal {D}^+$.
\end{lem}
\noindent {\bf Proof}: \,The proof in \cite{alsina2} only needs is
to be supplemented by the new notation of PN spaces with respect
to the
PM spaces as follows:\\
For any positive integer $n$, let $\sigma_n=\sum_{i=1}^n b_i$ and
let $\sigma=\sum_{i=1}^{\infty}b_i$. Then, one has
$$\tau_T^n(\epsilon_{b_1},\epsilon_{b_2},\cdots,\epsilon_{b_n})\geq\epsilon_{\sigma_n}>
\epsilon_{\sigma},
$$
whence $\tau_T^{\infty}\epsilon_{b_i}\geq\epsilon_{\sigma}$.
Consequently, since $G_{p_i}^i=\nu_{p_i}^i m_i\geq\epsilon_{b_i}$
for all $i\in \mathbf{N}$, one has $G_{\bar{p}}=\tau_T^{\infty}
G_{p_i}^i\geq\tau_T^{\infty}\epsilon_{b_i}\geq\epsilon_{\sigma}$.
Thus $G_{\bar{p}}(x)=1$ for $x>\sigma$.\qquad $\Box$
\bigskip

According to the result in Lemma 4 about the non trivial limit of
the infinite $\tau_T$--product, one has that $G_{\alpha \bar
{p}}$ and $G_{(1-\alpha)\bar{p}}$ are in $\mathcal {D}^+$ for
every $\alpha \in [0,1]$.
\bigskip

\noindent {\bf Example 4}. (A particular countable, but finite,
$\tau$--product) Let $(V,G)$ be the product of the Definition 8
for $i=1,2$, then $(V_1 \times V_2,\tau_T(\nu_{p_1}^1
m_1,\nu_{p_2}^2 m_2),\tau_T,\bf{T})$ is a PN
space.\\
To prove this it suffices to check axiom (iv): By Theorem 9 and
since ${\bf{T}}\gg\tau_T$ for every $t$--norm $T$ (see
\cite{schw-sk1}) one has
\begin{eqnarray*}
& G_{\bar {p}}=
\tau_T(\nu_{m_1}^1,\nu_{m_2}^2)\leq\tau_T({\bf{T}}(\nu_{\alpha
p_1}^1 m_1,\nu_{(1-\alpha)p_1}^1 m_1) ,{\bf{T}}(\nu_{\alpha
p_2}^2m_2,\nu_{(1-\alpha)p_2}^2m_2))\leq \\&
{\bf{T}}(\tau_T(\nu_{\alpha p_1}^1m_1,\nu_{\alpha
p_2}^2m_2)),\tau_T(\nu_{(1-\alpha)p_1}^1m_1,\nu_{(1-\alpha)p_2}^2m_2))
={\bf{T}}(G_{\alpha \bar{p}},G_{(1-\alpha)\bar{p}}) .\qquad \Box
\end{eqnarray*}
\bigskip

\noindent DEFINITION 11. Let $\{(V_i,\nu^i m_i,\bf{T},\bf{T})\}$
be a countable family of proper PN spaces. Then the
{\bf{T}}--product of this family is the pair $(V,G)$, where $V$
and $G$ are the same as in previous definition.
\bigskip

Just now one has the following question: Is the product $(V, G)$
a PN space under $\tau_T$ and $\tau_{T^*}$?.\\
\noindent First of all, can any member of the family $(V_i,\nu^i
m_i,\tau_T,\tau_{T^*})$ be a PN space?. We answer that question
in the negative because of the axiom (iv). Axiom (iii) works with
$m$ superadditive, but axiom (iv) needs $m$ to be subadditive.
The same function $m$ must to appear in second and third terms of
the following chain of inequalities

$$\nu_pm\leq \tau_{T^*}(\nu_{\alpha p},\nu_{(1-\alpha)p})m\leq \tau_{T^*}(\nu_{\alpha p}m,
\nu_{(1-\alpha)p}m),$$ and it would be meaningful that $m$ is
$\tau_{T^*}$--subadditive, whatever convergence factor one
introduces in $\tau_{T^*}$, which is absurd.\\
 In order to make the
reader's task easier the following result in PN spaces will be
needed, for the last section of this paper.
\bigskip

\begin{thm}\hspace{-0.22cm}\textbf{.} Let $\{(V_i,\nu^i m_i,\bf{T},\bf{T})\}$ be
a countable family of proper PN spaces. Then the pair $(V,G)$ is
a PN space under {\bf{T}}.
\end{thm}
\noindent {\bf Proof}:\, Clearly, for every $\bar{p},\bar{q}$ in
$V$, $G_{\bar{p}}=\epsilon_0$ if, and only if $\bar{p}=\theta$;
and $G_{\bar{p}}=G_{-\bar{p}}$.
\bigskip

\noindent (N3) By Lemma 4

\begin{eqnarray*}
 &  G_{\bar {p}+\bar{q}}= {\bf{T}}^{\infty}G_{p_i + q_i}\geq
{\bf{T}}^{\infty}({\bf{T}} (G_{p_i},G_{q_i}))=\\&
{\bf{T}}({\bf{T}}^{\infty}G_{p_i}^i,{\bf{T}}^{\infty}G_{q_i}^i)=
{\bf{T}}(G_{\bar{p}},G_{\bar{q}})
\end{eqnarray*}

\noindent (N4) By Lemma 4
$$G_{\bar{p}}\geq \epsilon_{\sigma}$$
for every $\bar{p}\in V$. Moreover,
\begin{eqnarray*}
 &  G_{\bar
 {p}}={\bf{T}}^{\infty}\nu_{p_i}^i m_i\leq{\bf{T}}^{\infty}({\bf{T}}(\nu_{\alpha p_i}^i,
\nu_{(1-\alpha)p_i}^i)\, m_i)=\\&
{\bf{T}}^{\infty}({\bf{T}}(\nu_{\alpha p_i}^i m_i,
\nu_{(1-\alpha)p_i}^i m_i)={\bf{T}}({\bf{T}}^{\infty}\nu_{\alpha
p_i}^i m_i, {\bf{T}}^{\infty}\nu_{(1-\alpha)p_i}^i
m_i)={\bf{T}}(G_{\alpha\bar{p}},G_{(1-\alpha)\bar{p}}).\qquad \Box
\end{eqnarray*}
\bigskip

\noindent {\bf Example 5}. Every member of the family
$\{(V_i,\nu_i,{\bf{\Pi}},{\bf{\Pi}})\}$, with the map $\nu_{p_i}$
defined via $ \nu_{p_i}:=exp\, (-\parallel p_i \parallel)$  is a
PN space as it is easy to check. Then, the family $\{(V_i,\nu_i
m_i,{\bf{\Pi}},{\bf{\Pi}})\}$ of $m_i$--transforms of the members
 of the above family is also a family of PN spaces, and his ${\bf{\Pi}}$--product
$(V,G)$ is a PN space under $\bf{\Pi}$, whence there exist
countable infinite products of PN spaces.
\bigskip

 \noindent DEFINITION 12. Let $\{(V_i, \nu^i,
\tau_i, \tau_i^*)\mid i\in \mathbf{N}\}$ be a countable family of
PN spaces. The $\Sigma$--product of this family is the pair $
\left(\prod_{i=1}^{\infty}V_i  , \,\nu^\Sigma\right)$ where
$\nu^\Sigma: \prod_{i=1}^{\infty}V_i \rightarrow {\Delta}^+$ is a
map given by

$$
\nu_{\bar {p}}^{\Sigma}:=\sum_{i=1}^{\infty}2^{-i}\nu_{p_i}^i,
$$
for every sequence $(p_i)=p\in \prod_{i=1}^{\infty}V_i$.\\

\noindent In order to simplify the notation we replace henceforth
$\nu_{\bar {p}}^{\Sigma}$ by $\nu_{\bar {p}}$.
\bigskip

\begin{thm}\hspace{-0.22cm}\textbf{.} Let $\{(V_i, \nu^i,\tau_i, \tau_i^*)\mid \in \mathbf{N}\}$ be
a countable family of PN spaces and let $\tau_i\geq\tau_W$ and
$\tau_i^*\leq \tau_W^*$ for all $i\in \mathbf{N}$, then the
$\Sigma$--product of this family denoted by
$$
\left(\prod_{i=1}^{\infty}V_i, \nu^\Sigma, \tau_W,
\tau_{W^*}\right)
$$
is a Menger space under $W$.
\end{thm}
\noindent The proof is similar to the one in Alsina [2]. It
suffices to apply Lemma 1 and Lemma 2.

\bigskip
\section {\bf Product topology for countable $\tau$--products}

 \bigskip
 In this section we want to show that the product topology and the
 strong topology in countable infinite $\tau$--products are not
 equal. Let us recall that this is what happened with the same type of
 $\tau$--products of PM spaces.

 \begin{thm}\hspace{-0.22cm}\textbf{.} Let each of the PN spaces
 $(V_i,\nu^i,\tau_i,\tau_i^*)$\, be endowed with the strong topology
 corresponding to $\nu^i, i\in \mathbf {N}$, and $\Delta^+$ with the
 topology of weak convergence. Then the product topology is weaker
 than the strong topology in $(V,G)$.
 \end{thm}
 \noindent
 {\bf Proof}: Let $U=\prod_{i=1}^{i=n} N_{p_i}(\epsilon_i)\times \prod_{i=n+1}^{\infty}V_i$
 be a standard neighborhood in the product topology. Choose $\epsilon=\min \{\epsilon_1,
 \epsilon_2,\ldots,\epsilon_n\}$ and let $\bar {q}\in N_{\bar
 {p}}(\epsilon)$. Then, since $G_{\bar {p}\,\bar
 {q}}\leq\nu_{q_i-p_i}^i$ for all $i\in \mathbf{N}$, one has

 $$
 1-\epsilon_i\leq 1-\epsilon<G_{\bar {p}\,\bar
 {q}}(\epsilon)\leq\nu_{q_i-p_i}^i(\epsilon)\leq \nu_{q_i-p_i}^i(\epsilon_i),
 $$

  \bigskip
  \noindent
  whence

 $$
 N_{\bar {p}}(\epsilon)\subset\prod_{i=1}^{\infty}N_{p_i}(\epsilon)\subsetneq U.
 \qquad \Box
 $$

 \bigskip
 In general, the two topologies are not equal. For, if this were the
 case, given $N_{\bar {p}}(\epsilon)$ there would exist a product
 neighborhood $U=\prod_{i=1}^{i=m} N_{p_i}(\epsilon_i)\times \prod_{i=m+1}^{\infty}V_i$
 such that $U\subset N_{\bar
 {p}}(\epsilon)\subset\prod_{i=1}^{\infty}N_{p_i}(\epsilon)$,
 which implies that $V_i=N_{p_i}(\epsilon)$ for all $i>m$, a very
 strong condition.
\bigskip

 \section {\bf Product topology  for $\Sigma$--products}

 \bigskip
 Contrary to what happens with countable infinite
 $\tau$--products,\, in $\Sigma$--products the two topologies
 are equal.\\

 \bigskip
 Let us recall that if $(V,\nu,\tau_W)$ is a PM space with $\tau_M$
 uniformly continuous, then for the strong neighborhoods $N_p(t),\, p\in V,
  \, t>0$,\, the following statements hold:
\begin{itemize}
   \item   If $q\in N_p(p)$, there exists a $t'>0$  such that $N_q(t')\subset N_p(t)$

   \item   If $p \ne q$, there exists a $t>0$  such that $N_p(t)\cap
   N_q(t)=\emptyset$
\end{itemize}

\bigskip
The proof of the following theorem is similar to that of Theorem
1.4 in \cite{alsina3}: only small changes in the notations are
needed.
\bigskip

 \begin{thm}\hspace{-0.22cm}\textbf{.} Let $\{(V_i,\nu^i,\tau_i,\tau_{W^*})|i\in \mathbf{N}\}$ and
 $V, \nu^\Sigma$ be as in theorem 8. Let each $V_i$ be endowed with
 the strong topology induced by $\nu^i$. Then the strong topology on
 $V$ induced by $\nu^\Sigma$ is the product topology.
 \end{thm}
 \noindent

The reason for the difference between Theorem 11 and 12 is easily
understood if one pays attention to the probabilistic
interpretation
of the $\epsilon$--neighborhood in the respective products spaces:\\
 \indent
If $N_{\bar {p}}(\epsilon)$ is a neighborhood in the
$\tau$--product, and $\bar {q}\in N_{\bar {p}}$ then, with
probability greater than $1-\epsilon$, all the components $p_i$
of $p$ are at a distance (the one associated to the norm in
$V_i$) less than $\epsilon$ from the corresponding $q_i$. On the
other hand, if $N_{\bar {p}}(\epsilon)$ is a neighborhood in the
$\Sigma$--product then $\bar {q}\in N_{\bar {p}}$ implies that,
with probability greater than $1-\epsilon$, at least one of the
components $p_i$ of $p$ is at a distance less than $\epsilon$
from the corresponding $q_i$.
\bigskip

\noindent
 {\bf Acknowledgments}\\

 The author wishes to thank Professors B. Schweizer and C. Sempi for
 their
interesting suggestions and comments as well as his patience while
reading the previous versions of this paper. The research of the
author was supported by grants from the Spanish
C.I.C.Y.T.(PB98-1010), and the Junta de Andaluc\'{\i}a.

\bigskip
\bigskip
\bigskip
\bigskip
\bigskip
\bigskip
\noindent{\it B. Lafuerza Guill\'{e}n}\\ \noindent Departamento de
Estad\'{\i}stica y Matem\'{a}tica Aplicada.\\Universidad de
Almer\'{\i}a, Spain\\ \noindent{\it e-mail}\, blafuerz@ual.es
\end{document}